\documentclass[11pt]{article}

\DeclareSymbolFont{rsfscript}{OMS}{rsfs}{m}{n}
\DeclareSymbolFontAlphabet{\mathrsfs}{rsfscript}
\renewcommand{\mathcal}{\mathrsfs}

\def\emline#1#2#3#4#5#6{%
       \put(#1,#2){\special{em:moveto}}%
       \put(#4,#5){\special{em:lineto}}}

\usepackage{latexsym}
\usepackage{amsmath}
\usepackage{amssymb}
\usepackage{amsfonts}
\usepackage{amsthm}

\theoremstyle{plain}

\newtheorem{theo}{Theorem}[section]
\newtheorem{lemma}[theo]{Lemma}
\newtheorem{cor}[theo]{Corollary}
\newtheorem{prop}[theo]{Proposition}

\newtheorem{rem}[theo]{Remark}
\newtheorem*{remo}{Remark}
\newtheorem{exams}[theo]{Examples}

\DeclareMathOperator{\loc}{loc}
\DeclareMathOperator{\rind}{ind}
\DeclareMathOperator{\spa}{span}

\DeclareMathOperator{\id}{\rm id}
\DeclareMathOperator{\proj}{\rm proj}
\DeclareMathOperator{\vol}{\rm vol}

\newcommand{\Proof}{\par\bigskip\noindent{\sc Proof: }}
\newcommand{\eop}{\nopagebreak\hspace*{\fill}$\blacksquare$\bigskip}

\newcommand{\C}{{\mathbb C}}

\newcommand{\R}{{\mathbb R}}
\newcommand{\N}{{\mathbb N}}
\newcommand{\Z}{{\mathbb Z}}

\newcommand{\noi}{\noindent}
\newcommand{\bn}{\bigskip\noindent}

\newcommand{\sn}{\smallskip\noindent}
\newcommand{\bsk}{\bigskip}

\newcommand{\ssk}{\smallskip}
\newcommand{\q}{\quad}
\newcommand{\e}{\ell}
\newcommand{\wt}{\widetilde}
\newcommand{\wh}{\widehat}
\newcommand{\ov}{\overline}
\newcommand{\os}{\overset}
\newcommand{\lora}{\longrightarrow}
\newcommand{\hra}{\hookrightarrow}
\newcommand{\car}{\curvearrowright}

\newcommand{\ma}{\mathop}
\newcommand{\li}{\limits_}
\newcommand{\tsum}{{\textstyle{\sum}}}

\newcommand{\suk}{\tsum(u_k)}
\newcommand{\dvd}{\| \cdot \|}
\newcommand{\ce}{\mathcal E}
\newcommand{\cd}{\mathcal D}
\newcommand{\cb}{\mathcal B}

\newcommand{\cl}{\mathcal L}

\newcommand{\fv}{\mathfrak V}
\newcommand{\rn}{\R^n}
\newcommand{\vn}{v_n}
\newcommand{\varsub}{\varsubsetneqq}
\newcommand{\skinf}{\sum_{k=1}^\infty}
\newcommand{\es}{\emptyset}
\newcommand{\ex}{\exists}
\newcommand{\ssp}{\supset}
\newcommand{\fa}{\forall}
\newcommand{\Om}{\Omega}

\newcommand{\la}{\langle}
\newcommand{\ra}{\rangle}
\newcommand{\bl}{\bigl(}
\newcommand{\br}{\bigr)}
\newcommand{\Bl}{\Bigl(}
\newcommand{\Br}{\Bigr)}
\newcommand{\bv}{\big|}

\newcommand{\bmv}{\bigm|}
\newcommand{\Bmv}{\Bigm|}

\newcommand{\blk}{\bigl\{}
\newcommand{\brk}{\bigr\}}
\newcommand{\Blk}{\Bigl\{}
\newcommand{\Brk}{\Bigr\}}

\newcommand{\Ban}{Banaszczyk}
\newcommand{\Fre}{Fr\'echet}

\newcommand{\ci}{\scriptstyle\circ}
\title{\bf The Levy-Steinitz rearrangement theorem\\[1ex] for
duals of metrizable spaces}
\author{%
{\small Jos\'e Bonet}\footnote{The research of the first author was
partially supported by DGES Project PB97-0333, and the second one by a
grant of the Ministerio de Educati\'on y Cultura of Spain (Ref:
Sab 1995-0736)}\\
{\small Departamento de Matem\'atica Aplicada}\\
{\small Universidad Polit\'ecnica de Valencia}\\
{\small E-46022 Valencia, Spain}\\
{\small E-mail: jbonet@pleiades.upv.es}\\[2ex]
{\small Andreas Defant}\\
{\small Fachbereich Mathematik}\\
{\small Universit\"at Oldenburg}\\
{\small D-26111 Oldenburg, Germany}\\
{\small E-mail: defant@mathematik.uni-oldenburg.de}
}
\date{}

\begin{document}
\maketitle

\begin{abstract}
Extending the Levy-Steinitz rearrangement theorem in $\R^n$, which in turn
extended Riemann's theorem, Banaszczyk proved in 1990/93 that a metrizable,
locally convex space is nuclear if and only if the domain of sums of every
convergent series (i.e. the set of all elements in the space which are sums
of a convergent rearrangement of the series) is a translate of a closed
subspace of a special form. In this paper we present an apparently complete
analysis of the domains of sums of convergent series in duals of metrizable
spaces or, more generally, in (DF)-spaces in the sense of Grothendieck.
\end{abstract}

\section*{Introduction}

For a convergent series $\tsum(u_k)$ in a locally convex space $E$ the
domain of sums $S\bl\tsum(u_k)\br$ is the set of all $x\in E$ which can be
obtained as the sum of a convergent rearrangement of $\tsum(u_k)$. In terms
of this notion Riemann's famous rearrangement theorem states that in the
real line $\R$ the domains of sums are either single points or
coincide with the whole line. Later Levy \cite{L} and Steinitz \cite{S}
extended Riemann's
result to finite dimensional spaces: domains of sums in $\R^n$ are affine
subspaces -- more precisely, for each convergent series $\tsum(u_k)$ in
$\R^n$
\[
  S\bl\tsum(u_k)\br = \sum^\infty_{k=1} u_k +\Blk y\in \R^n \Bmv
  \la x,y\ra = 0 \text{ for all } x\in \R^n \text{ with }
  \sum^{\infty}_{k=1} \bv \la x,u_k\ra \bv < \infty \Brk \, ;
\]
we refer to this result as the ``Levy-Steinitz theorem''
(see \cite{Ro} and \cite{KK} for a proof). For the
``state of art'' of this theorem
in infinite dimensions (in particular, in Banach spaces)
see the recent monograph \cite{KK} of M.I.~Kadets and
V.M.~Kadets.

The following notation is basic to the understanding of domains of sums in
infinite
dimensional spaces: for each convergent series $\tsum(u_k)$ in a locally
convex space $E$ define the set
\begin{align*}
   \Gamma\bl \tsum(u_k)\br\;\:\, & := \blk x'\in E'\bmv \bl x'(u_k)\br \in
    \e_1\brk\subset E'\\
\intertext{and its polar}
   \Gamma \bl \tsum(u_k)\br^{\bot}   & :=  \blk x\in E \bmv x'(x) = 0
   \text{ for all } x'\in \Gamma (\tsum (u_k)\br \brk \subset E\, .
\end{align*}
Obviously, both sets are subspaces of $E$ and $E'$, respectively,
and the second one is even closed. It
is an easy exercise to check the following remark which is useful later
(for $E = \R^n$ see also \cite[Thm.1]{CC}):

\begin{remo}
Let $E$ be a locally convex space and $\tsum(u_k)$ a convergent series.
Then for each $x\in E$ the following are equivalent:
\begin{itemize}
\item[{\rm (1)}]
 $x\in \ma{\sum}\li{k=1}^{\infty}u_k + \Gamma \bl \tsum
(u_k)\br^{\bot}$
\item[{\rm (2)}]
 $\forall x'\in E'\; \ex \text{ permutation } \pi \text{ of } \N: x'(x) =
\ma{\sum}\li {k=1}^\infty x'(u_{\pi(k)})\, $.
\end{itemize}
\end{remo}
In two papers \cite{B1}, \cite{B2} from 1990 and 1993 \Ban \ proved the
following
extension of the Levy-Steinitz theorem -- a result which here will be
quoted as ``\Ban 's rearrangement theorem'':
\begin{quote}
\it
A metrizable, locally convex space $E$ is nuclear if and only if the domain
of sums of each convergent series $\tsum(u_k)$ in $E$ is given by the
formula
\[
   S\bl\tsum(u_k)\br = \sum^\infty_{k=1}u_k +\, \Gamma\bl\suk\br^{\bot}\, ;
\]
in particular, $S\bl\suk\br$ is a closed, affine subspace of $E$, and for
each $x\in E$:
\[
   x\in S\bl\suk\br \text{ iff } \forall x'\in E': x'(x)\in
   S\Bl\tsum\bl x'(u_k)\br\Br\, .
\]
\end{quote}

The aim of this paper is to investigate the domains of sums in duals of
metrizable spaces or, more generally, in (DF)-spaces in the sense of
Grothendieck (in fact,
in \cite[p.101]{KK} the authors write that ``Steinitz-type problems for
locally convex spaces are poorly investigated'').
Almost all natural spaces of functions and distributions are metrizable
spaces, or duals of metrizable spaces or, if not of this type, then at
least ``generated'' by spaces of this type.
The only infinite dimensional (DF)-space for which
the Levy-Steinitz theorem was known to hold is the countable direct sum
$\varphi:=\bigoplus_{\N}\R$ of copies of the scalar field $\R$
(this space
is the strong dual of the nuclear \Fre \ space $w:= \prod_{\N}\R$);
we refer to \cite[Ex. 8.3.3]{KK}.
Implementing the basic ingredients of \cite{B1} and \cite{B2} in an
alternative way into locally convex spaces, we are able to say
precisely which parts of \Ban 's theorem transfer to duals of metrizable
spaces and which don't.

Let us briefly describe the content of this article. After a short sketch
of
\Ban 's theorem in section 1, we show in section 2 the following
rearrangement theorem: if
$E$ is the dual of a nuclear, metrizable space, then each domain of
sums $S\bl\suk\br$ is the translate of a subspace of $E$ of a special
form:
\[
  S\bl\suk\br = \sum_{k=1}^{\infty}u_k + \, \Gamma_{\loc}^{\bot}
  \bl\suk\br\, .
\]
If $E$ is not isomorphic to $\varphi$, then there is a convergent series
whose domain of sums is not closed -- in particular, the Levy-Steinitz
theorem in its original formulation does not hold in the dual of a nuclear,
metrizable space unless this space equals the trivial space $\varphi$. In
section 3 it
is proved conversely that if $E$ is a dual of a nonnuclear metrizable
space, then there is a convergent series whose domain of sums does not
have this special form.

An analysis of our approach using local convergence and bounded sets
permits to show that in a much larger class of (nonmetrizable) spaces the
domain of sums of all convergent series are affine subspaces;
this class includes for example: the space $H(K)$ of germs of holomorphic
functions on a
compact set K, the space of Schwartz tempered distributions $S'$, the space
$\ce'(\Om)$ of distributions with compact support, the space $\cd(\Om)$ of
test functions for distributions,
the space $\cd'(\Om)$ of distributions, the space $A(\Om)$ of real analytic
functions
and its dual $A'(\Om)$, for an open subset $\Om$ of $\R^N$.

\section*{Preliminaries}

We shall use standard notation and notions from the theory of locally
convex spaces and Banach spaces (see e.g. \cite{BPC}, \cite{G}, \cite{Ja},
\cite{Ju}
and \cite{MV}); for the theories of operator ideals and $s$-numbers see
\cite{K},
\cite{P2} and \cite{P3}. All locally convex spaces $E$ are real, by $cs\, (E)$
we denote the collection of all continuous seminorms, and by $\cb(E)$ all
absolutely convex, closed and bounded sets of $E$.
The dual $E'$ endowed with the topology of uniform convergence on all
bounded sets is denoted by $E'_b$.
A sequence $(p_n)$ in $cs\, (E)$
(or $(B_n)$ in $\cb(E)$) is said to be a fundamental system
of $cs\, (E)$ (or $\cb(E)$) whenever it is increasing (natural
order) and for each $p\in cs\, (E)$ (or $B\in \cb(E)$) there is $n$ such
that $p\le p_n$ (or $B\subset B_n$). Recall that a locally convex space
is metrizable if and only if it allows a fundamental system $(p_n)$ of
$cs\, (E)$, and a locally convex space $E$ is said to be a (DF)-space if it
has a fundamental system $(B_n)$ of $\cb(E)$ and, moreover, each
intersection of any sequence of closed, absolutely convex
zero neighbourhoods in $E$ is again a zero neighbourhood, provided it
absorbs all bounded sets. Strong duals of metrizable spaces are (DF).
For a seminorm $p$ on  a vector space $E$ we write $E_p$ for the normed
space $(E/\ker p,\|\cdot\|_p)$, $\|x+\ker p\|_p := p(x)$, and for an
absolutely convex set $B\subset E$ the seminormed space $(\spa
B,\, m_B)$, $m_B$ the Minkowski gauge functional, is denoted by $E_B$.
The canonical surjection $E\lora E_p$ and injection $E_B\hra E$ are denoted
by $\pi_p$ and $i_B$, respectively. Clearly, we have $\os{\ci}{B}_{E_B}
\subset B\subset B_{E_B}$ for the open and closed unit ball of $E_B$. The
set $B$ is said to be a Banach or Hilbert disc whenever $E_B$ is a Banach
or Hilbert space. For a locally convex space $E$ and $B\in \cb(E)$ the
natural embedding $i_B$ is continuous, hence $B = B_{E_B}$ ($B$ is closed
in $E_B$). In complete spaces $E$ all $B\in \cb(E)$ are Banach discs
(see \cite[5.1.27]{BPC} or \cite[23.14]{MV}).
A locally convex
space $E$ is called nuclear if for each $p\in cs\, (E)$ there is $p\le
q\in cs\, (E)$ such that the canonical linking map $\pi_q^p: E_q\lora E_p$
is nuclear, and $E$ is said to be co-nuclear if its dual is nuclear or,
equivalently, for each $B\in \cb(E)$ there is $B\subset C\in \cb(E)$ such
that the embedding $i^C_B: E_B\hra E_C$ is nuclear. Metrizable spaces and
(DF)-spaces are nuclear if and only if they are co-nuclear, and nuclear
\Fre \ spaces and complete nuclear (DF)-spaces are in a one-to-one relation
with respect to the correspondence $E\lora E'_b$ (see \cite[7.8.2]{Ju}).
For further
details on nuclearity we refer to \cite{Ju}, \cite{MV} or \cite{P1}.
A sequence $(x_n)$ in a locally convex space $E$ is said to converge locally
whenever it is contained and converges in some $E_B$, $B\in \cb(E)$.
A locally convex space $E$ satisfies the strict Mackey condition
(see e.g. \cite[5.1.29]{BPC}) if given any $B\in\cb (E)$ there is $C\in \cb
(E)$,
$C\supset B$, such that $E$ and $E_C$ induce the same topology on $B$.
Clearly, in every locally convex space with the strict Mackey condition
every convergent sequence converges locally. Every metrizable space
satisfies the strict Mackey condition (\cite[5.1.30]{BPC}) as well as every
complete nuclear (DF)-space (use e.g. \cite[7.3.7]{Ju}).

%1
\section{On \Ban 's rearrangement theorem for nuclear\\ \Fre \ spaces}

For later use (and in order to set up some further notation) we want to
give a short sketch of the proof of \Ban 's rearrangement theorem (which
still follows the line of the original proof of Steinitz for $\R^n$).

\bsk
Let $\tsum(u_k)$ be a convergent series in a locally convex space $E$. The
extended domain of sums $S^e\bl\suk\br$ is defined as the set of all
$x\in E$ which appear as limits of convergent subsequences of
rearrangements $\tsum\bl u_{\pi(k)}\br$. Moreover, let
\begin{align*}
  A\bl \suk\br & := \bigcap_m Z_m\bl\suk\br\\
  Q\bl \suk\br & := \bigcap_m co \Bl Z_m\bl\suk\br\Br\, ,
\end{align*}
where
\[
  Z_m\bl\suk\br := \Blk\sum_Iu_k \bmv I\subset \blk m,m+1,\ldots \brk
                   \text{ finite }\Brk\, ,
\]
and
\[
  A_E\bl \suk \br := \bigcap_m \ov{Z}_m^E \text{ and }
  Q_E\bl \suk \br := \bigcap_m \ov{co\, Z_m}^E\, .
\]
The proof of part (a) in the following result is elementary (for a more
general formulation
for abelian Hausdorff groups
see \cite[Ch.III, \S 5, Ex.3]{Bo}), and (b) is a consequence of
an easy Hahn-Banach argument (check the proof of \cite[Lemma 6]{B1} or
see \cite{Bl}).

\begin{lemma}\label{lemma1.1}
Let $\suk$ be a convergent series in a locally convex space $E$. Then
\begin{itemize}
\item[{\rm (a)}]
 $S^e\bl \suk\br = \ma{\sum}\li{k=1}^\infty u_k + A_E\bl\suk\br$,
    whenever $E$ is metrizable.
\item[{\rm (b)}]
 $Q_E\bl \suk\br = \Gamma^{\bot}\bl \suk\br$.
\end{itemize}
\end{lemma}

\sn
In the metrizable case this obviously gives the following chain of
equalities and inclusions:

\bsk
\begin{lemma}\label{lemma1.2}
Let $\tsum(u_k)$ be a convergent series in a metrizable space $E$. Then
\begin{align}
 S\bl\suk\br & \subset S^e\bl\suk\br
                  \tag*{(B1)}\\
             & = \sum_{k=1}^\infty u_k + A_E\bl\suk\br
                 \notag\\
             & \subset \sum_{k=1}^\infty u_k + Q_E\bl\suk\br
                 \tag*{(B2)}\\
             & = \sum_{k=1}^\infty u_k + \Gamma^{\bot}\bl \suk\br\, .
                 \notag
\end{align}
\end{lemma}

\bn
In nuclear metrizable spaces the inclusions (B1) and (B2)
are even equalities -- these are the crucial steps in \Ban 's paper
\cite{B1}. The proofs are consequences of the following two lemmas
(which are variants of \cite[Lemma 4 and Lemma 8]{B1}). Here $HS(T:
X\lora Y)$ means the Hilbert-Schmidt norm of an operator $T$ acting between
the Hilbert spaces $X$ and $Y$.

\begin{lemma}[the ``lemma of rounding off coefficients'']\label{lemma1.3}
Let $H_1$ and $H_2$ be two Hil\-bert spaces such that $H_1 \subset H_2$
and
\[
  HS\,(\id: H_1\hra H_2)\le 1\, .
\]
Then for $y_1,\ldots, y_s \in B_{H_1}$ and $y\in co\, \Blk \ma{\sum}\li J
y_k
\bmv J \subset \blk 1,\ldots, s\brk\Brk$ there is a finite
set\linebreak[4] $I\subset
\blk 1,\ldots, s\brk$ such that
\[
   \sum_I y_k-y \in B_{H_2}\, .
\]
\end{lemma}

\begin{lemma}[the ``permutation lemma'']\label{lemma1.4}
Let $H_k$, $k = 1,2,3$, be three Hilbert spaces such that $H_1\subset
H_2\subset H_3$ and
\begin{align*}
  HS\, (\id: H_1\hra H_2) &  \le 1\, ,\\
  HS\, (\id: H_2\hra H_3) & \le 1/2 \, .
\end{align*}
Then for $v_1,\ldots, v_s\in B_{H_1}$ and $a\in B_{H_2}$ with $a+
\ma{\sum}\li {k=1}^s v_k\in B_{H_2}$
there is a permutation $\sigma$ of $\{ 1,\ldots,s\}$ such that
\[
   a + \sum^m_{k=1} v_{\sigma (k)}\in B_{H_3} \text{ for all }
   m=1,\ldots,s\, .
\]
\end{lemma}

%2
\section{The rearrangement theorem for nuclear (DF)-spaces}

In this section we will analyze the domain of sums of convergent series in
complete nuclear (DF)-spaces $E$. Clearly, $\rn$ and $\varphi = \bigoplus_
{\N}\R$ are examples of complete nuclear (DF)-spaces for which the
Levy-Steinitz theorem
holds -- and we will show that these are in fact the only ones with this
property, although all domains of sums in complete nuclear (DF)-spaces turn
out to be affine subspaces.

\bsk
Take a convergent series $\suk$ in a nuclear (DF)-space $E$, and recall
from the preliminaries that
it converges locally, i.e., there is some $B\in \cb (E)$ such that $\suk$
converges in $E_B$. Moreover, let $(B_n)$ be a fundamental system of $\cb
(E)$ such that $B_1 = B$. Then it follows from Lemma \ref{lemma1.2} that,
for each $n$,
\begin{align*}
  S_{B_n}\bl\suk\br & \subset S_{B_n}^e\bl\suk\br\\
                    & = \sum_{k=1}^\infty u_k + A_{B_n}\bl\suk\br\\
                    & \subset \sum_{k=1}^\infty u_k + Q_{B_n}\bl\suk\br
                      = \sum_{k=1}^\infty u_k + \Gamma^{\bot}_{B_n}
                        \bl\suk\br\, ,
\end{align*}
where the index $B_n$ indicates that we consider all sets involved with
respect to the convergent series $\suk$ in $E_{B_n}$. We obtain as an
immediate consequence that
\begin{align}
 S\bl\suk\br & = \bigcup_n S_{B_n}\bl\suk\br \notag\\
             & \subset \bigcup_n S_{B_n}^e\bl\suk\br
                  \tag*{(B1')}\\
             & = \sum_{k=1}^\infty u_k + \bigcup_n A_{B_n}\bl\suk\br
                  \notag\\
             & \subset \sum_{k=1}^\infty u_k + \bigcup_n Q_{B_n}
                  \bl\suk\br
                 \tag*{(B2')}\\
             & = \sum_{k=1}^\infty u_k + \bigcup_n \Gamma_{B_n}^{\bot}
                  \bl\suk\br\, .
                  \notag
\end{align}
Define
\[
  \Gamma_{\loc}^{\bot}\bl\suk\br := \bigcup_n
  \Gamma_{B_n}^{\bot}\bl\suk\br
\]
and observe that this set is independent of the bounded set $B$
and the fundamental system $(B_n)$ we chose at the beginning. Since the
sets $\Gamma^{\bot}_{B_n}$ form an increasing family of subspaces of $E$,
the set $\Gamma_{\loc}^{\bot}$ is even a subspace of $E$.

\begin{theo}\label{theo2.1}
Let $E$ be a complete nuclear (DF)-space.
\begin{itemize}
\item[{\rm (a)}]
 For each convergent series $\suk$ in $E$ the domain of sums $S\bl\suk\br$
is an affine subspace of $E$; more precisely:
\[
   S\bl\suk\br = \sum_{k=1}^\infty u_k + \Gamma_{\loc}^{\bot}\,
                 \bl\suk\br\, .
\]
\item[{\rm (b)}]
 Assume $E$ to be infinite dimensional and not isomorphic to $\varphi$.
Then
there is a convergent series $\suk$ in $E$ such that the domain of sums
$S\bl \suk\br$ is a nonclosed subspace of $E$; in particular:
\[
   S\bl\suk\br = \sum_{k=1}^\infty u_k + \Gamma^{\bot}_{\loc}\, \bl\suk\br
   \varsub \skinf u_k + \Gamma^{\bot}\bl\suk\br\, .
\]
\end{itemize}
\end{theo}
As in the metrizable case statement (a) will follow from the fact that the
above inclusions (B1') and (B2') are even equalities. To prove this,
it is necessary to adapt the lemma of rounding off coefficients and the
permutation lemma (1.3 and 1.4 of section 1) to the
(DF)-setting, and this is done in the following two lemmas:

\begin{lemma}\label{lemma2.2}
Let $E$ be a vector space and $B_0\subset B\subset C$ three Hilbert discs
in $E$ such that
\begin{align*}
   & HS\, \bl i^B_{B_0}: E_{B_0}\hra E_B \br \le 1\, ,\\
   & HS\, \bl i^C_B: E_B\hra E_C \br \le 1/2\, .
\end{align*}
Then for each convergent series $\suk$ in $E_{B_0}$
\[
   S_B^e\bl\suk\br \subset S_C\bl\suk\br\, .
\]
\end{lemma}

\Proof
Take
\[
   x = \lim_{m\to \infty} \sum_{k=1}^{j_m} u_{\pi(k)} \in
       S^e_B\bl \suk \br
\]
(convergence in $E_B$) and extract an increasing subsequence
$(j_{m(\e)})_{\e\ge
2}$ of
$(j_m)$ such that for all $\e\ge 2$
\begin{align*}
 &   x_{\pi(k)}\in 1/\e\, B_0 \text{ for all } k= j_{m(\ell)}+1,\ldots,
         j_{m(\ell+1)}\, ,\\
 &   x - \sum_{k=1}^{j_{m(\e)}} u_{\pi(k)}\in 1/\e \,B\, ,\\
 &   x - \sum_{k=1}^{j_{m(\ell+1)}} u_{\pi(k)}\in 1/\e \,B \, ;
\end{align*}
this is possible: take first $k_0$ such that
\[
  u_{\pi(k)}\in 1/2\, B_0 \text{ for all } k\ge k_0\, ,
\]
and then $m_0$ such that $j_{m_0} \ge k_0$ and
\[
  x - \sum^{j_m}_{k=1} u_{\pi(k)}\in 1/2\, B \text{ for all $m\ge m_0$.}
\]
Put $m(2):= m_0$. Then for $m\ge m(2)$
\begin{align*}
  &  u_{\pi(k)}\in 1/2\, B_0 \text{ for all } k\ge j_{m(2)}\, ,\\
  &  x - \sum^{j_m}_{k=1} u_{\pi(k)}\in 1/2\, B\text{ for all } m\ge m(2)\,
.
\end{align*}
Select now $m(3) > m(2)$ such that
\begin{align*}
   & u_{\pi(k)}\in 1/3 \, B_0 \text{ for all } k\ge j_{m(3)}\, ,\\
   & x - \sum^{j_m}_{k=1} u_{\pi(k)}\in 1/3\, B\text{ for all } m\ge m(3)\,
.
\end{align*}
Then
\begin{align*}
   & u_{\pi(k)}\in 1/2 \, B_0 \text{ for all } k =
                 j_{m(2)}+1,\ldots,j_{m(3)}\, ,\\
   & x - \sum^{j_{m(2)}}_{k=1} u_{\pi(k)}\in 1/2\, B\, ,\\
   & x - \sum^{j_{m(3)}}_{k=1} u_{\pi(k)}\in 1/2\, B\, ,\\
\end{align*}
etc. $\ldots\;$. Consider now the Hilbert discs
\[
  1/\e \, B_0 \subset 1/\e\, B \subset 1/\e\, C\, ,
\]
Clearly,
\[
  HS\bl E_{1/\e\, B_0} \hra E_{1/\e\,B} \br \le 1 \q \text{and} \q
  HS\bl E_{1/\e\, B} \hra E_{1/\e\,C}\br \le 1/2 \, .
\]
Hence, by Lemma \ref{lemma1.4}, for each $\ell \ge 2$ there exists a
permutation $\sigma_{\e}$ of
$\pi\blk j_{m(\e)} +1,\ldots, j_{m(\e+1)}\brk$ such that
 for all $m = j_{m(\e)}+1,\ldots, j_{m(\e+1)}$
\[
   x - \sum^{j_{m(\e)}}_{k=1} u_{\pi(k)} -
   \sum^m_{k=j_{m(\e)}+1} u_{\sigma_{\e}\bl\pi(k)\br} \in 1/\e\, C\, .
\]
Clearly, this gives a permutation $\varrho$ of $\N$ such that (convergence
in $E_C$)
\[
   x = \skinf u_{\varrho(k)}\, ,
\]
hence $x\in S_C\bl \suk\br\, .$
\eop

\begin{lemma}\label{lemma2.3}
Let $E$ be a vector space and $B_0\subset B$ two Hilbert discs in $E$ such
that
\[
   HS \Bl i^B_{B_0}: E_{B_0} \hra E_B\Br \le 1 \, .
\]
Then for each convergent series $\suk$ in $E_{B_0}$
\[
   A_B \bl\suk\br = Q_B\bl\suk\br\, .
\]
\end{lemma}

\Proof
According to the definitions $A_B\subset Q_B$. Take $x\not\in A_B$.
This means, there is $m_0$ such that
\[
  x\not\in \ov{Z}_m^{E_B} \text{ for all } m\ge m_0\, .
\]
Choose $\varepsilon > 0$ such that
\[
   (x+\varepsilon B)\cap Z_{m_0} = \es \, ,
\]
and $m_1 > m_0$ such that
\[
    u_k\in \varepsilon/2\, B_0 \text{ for all } k \ge m_1\, .
\]
We will show that
\[
    (x + \varepsilon/2 \, B) \cap co\, Z_{m_1} = \es \,  ,
\]
hence $x\not\in Q_B$. If this is not the case, then there is
\[
   z\in (x + \varepsilon/2\, B)\cap co\, Z_{m_1}\, .
\]
Since $HS \bl E_{\varepsilon/2\, B_o} \hra  E_{\varepsilon/2\, B}\br
\le 1$, by Lemma \ref{lemma1.3} there is $\wt{z}\in Z_{m_1}$ with
\[
  z - \wt{z}\in \varepsilon/2\, B\, ,
\]
hence
\[
  \wt{z}\in \bl z+ \varepsilon/2\, B\br\cap Z_{m_1}\subset
   (x+\varepsilon\, B)\cap Z_{m_0}\, ,
\]
%\nopagebreak
a contradiction.
\eop

\bsk
Now we are prepared to give the {\sc proof} of part (a) of Theorem
\ref{theo2.1}: Let $\suk$
be a convergent series in a complete nuclear (DF)-space $E$. Take a
fundamental
system $(B_n)$ of $\cb (E)$ such that $\suk$ converges in $E_{B_1}$, all
$B_n$ are Hilbert discs and
\[
   HS \Bl i^{B_{n+1}}_{B_n}: E_{B_n}\hra E_{B_n+1}\Br\le 1/2 \text{ for
          all } n
\]
(see the preliminaries and combine e.g. \cite[7.8.2(4)]{Ju} with
\cite[7.6.3(3)]{Ju}).
Then the inclusion (B1') by Lemma \ref{lemma2.2} and the inclusion
(B2') by Lemma \ref{lemma2.3} are equalities which gives the claim.
\hfill$\blacksquare$

\bsk
This proof obviously allows collecting a bit more of information about
domains of sums in nuclear (DF)-spaces:

\begin{rem}\label{rem2.4}
Let $\suk$ be a convergent series in a complete nuclear (DF)-space $E$, and
$(B_n)$
a fundamental system of $\cb(E)$ such that $\suk$ converges in $E_{B_1}$.
Then
\begin{align*}
  S\bl\suk\br & = \bigcup_n \, S_{B_n} = \bigcup_n\, S^e_{B_n}\\[1ex]
              & = \skinf u_k + \bigcup_n\, A_{B_n} = \skinf u_k + \bigcup_n
                   \, Q_{B_n}\\[1ex]
              & = \skinf u_k + \bigcup_n \, \Gamma^{\bot}_{B_n} = \skinf u_k
                  + \Gamma^{\bot}_{\loc}\bl\suk\br\, .
\end{align*}
\end{rem}

The proof is clear from what was said before (note that all unions are
independent of the fundamental system $(B_n)$ which hence can be chosen as
in the above proof of 2.1 (a)).

\bsk
It remains to give a {\sc Proof} of part (b) of Theorem \ref{theo2.1}: The
nuclear
\Fre\ space $F := E'_b$ is not isomorphic to $w := \prod_{\N}\R$ (recall
that $E$ is reflexive).
Accordingly, we can apply \cite[Theorem 5]{DV} to conclude the existence of
a quotient $G$ of $F$ which has a continuous norm $\| \cdot \|$, but is not
countably normable (a projective limit with injective connecting mappings).
We choose a fundamental system $(\|\cdot \|_n)$ of seminorms in $G$ such
that
$\dvd \le \dvd_1$ and all the completions of the normed spaces $(G,
\dvd_n)$ are separable
and reflexive (see e.g. \cite[7.6.3]{Ju}). Since $F$ and $G$ are nuclear
\Fre\
spaces, it follows that $H = G'_b$ is a closed, topological subspace of $E
= F'_b$. For each $n$ define
\[
  B_n := \{ x\in G\mid \| x \|_n\le 1\}^0 \subset H \text{ (polar in $H$)
    and } H_n := H_{B_n}\, .
\]
We have $H:= \rind_n\, H_n$, and since $\dvd_1$ is a norm on $G$, the Banach
space $H_1$ is dense in $H$. As $G$ is not countably normable, we can apply
\cite[2.9 and 2.6]{BMT} to conclude
\[
  \bigcup_n \ov{H}^{H_n}_1 \varsub \ov{H}^H_1 = H\, .
\]
By \cite[Ex. 3.1.5]{KK}, since $H_1$ is separable, there is a convergent
series $\suk$ in $H_1$ such that its domain of sums in $H_1$ coincides with
$H_1$. The domain of sums $S\bl\suk\br$ in $H$ is a subspace of $H$, which
contains $H_1$ but is not closed; indeed, if it is closed, it must coincide
with $H$. On the other hand,
\[
  S\bl\suk\br \subset \bigcup_n\, S_{B_n}\bl\suk\br \subset \bigcup_n
  \ov{H}_1^{H_n}\ne H\, .
\]
This is a contradiction.
\eop

In view of the remark made in the introduction the following consequence of
2.1 (b) seems to be notable:

\begin{cor}\label{cor2.5}
Let $E$ be a complete and infinite dimensional, nuclear (DF)-space which is
not isomorphic to $\varphi$. Then there is a convergent series $\suk$ in
$E$ and some $x\in E$ such that for all $x'\in E'$ there is a permutation
$\pi$ of $\N$ with
\[
  x'(x) = \skinf x'\bl u_{\pi(k)})\, ,
\]
but $x\not\in S\bl\suk\br$.
\end{cor}

We finish this section showing that our approach to the Levy-Steinitz type
theorems via local convergence and bounded sets covers a much larger class
of spaces than only nuclear (DF)-spaces -- it turns out that
co-nuclearity and not nuclearity is the appropriate assumption needed.

\bsk
Let $E$ be a locally convex space in which every sequence converges
locally, and let $\suk$ be a series in $E$ which converges in,
say, $E_{B_0}$.
Denote all absolutely convex and bounded supersets of $B_0$ by $\cb$.
Clearly,
\pagebreak[4]
\begin{align*}
 S\bl\suk\br & = \bigcup_{\cb}\, S_B \subset \bigcup_{\cb}\,
                 S^e_B\\[1ex]
             & = \skinf u_k + \bigcup_{\cb}\, A_B\\[1ex]
             & \subset \skinf u_k +\bigcup_{\cb} \, Q_B\\[1ex]
             & = \skinf u_k + \bigcup_{\cb} \Gamma^{\bot}_B =:
                 \skinf u_k + \Gamma^{\bot}_{\loc}\bl\suk\br\, ,
\end{align*}
where the notation of all new symbols is obvious; again the set
$\Gamma^{\bot}_{\loc}$ is independent of the choice of $B_0$, and it can
easily be shown that it is even a subspace of $E$. Another application of
2.2 and 2.3 gives the following extension of part (a) of Theorem
\ref{theo2.1}:

\begin{theo}\label{theo2.6}
Let $E$ be a complete, co-nuclear space for which each convergent sequence
converges locally. Then the domain of sums of each convergent series $\suk$
in $E$ is an affine subspace,
\[
   S\bl\suk\br = \skinf u_k + \Gamma^{\bot}_{\loc}\bl\suk\br\, .
\]
\end{theo}

Before we give a list of new examples of spaces such that each domain of
sums is affine, let us show that the preceding result even covers \Ban 's
rearrangement theorem for metrizable spaces. For this it suffices to check
the following remark, since a metrizable space is nuclear if and only if it
is co-nuclear, and each convergent sequence in a metrizable space converges
locally (see the preliminaries).

\bsk
\begin{rem}\label{rem2.7}
For each convergent series $\suk$ in a metrizable, locally convex space $E$
\[
  \Gamma^{\bot}\bl\suk\br = \Gamma^{\bot}_{\loc}\bl\suk\br \, .
\]
\end{rem}

\Proof
Clearly
\begin{align*}
\Gamma^{\bot}_{\loc}  & = \bigcup_{\cb} \Gamma^{\bot}_B =
                        \bigcup_{\cb}\,\bigcap_m \ov{co\, Z_m}^{E_B}\\[1ex]
                     & \subset \bigcap_m \ov{co\, Z_m}^E =
                        \Gamma^{\bot}\, ,
\end{align*}
hence take $x\in \ma{\bigcap}\li m \ov{co\, Z_m}^E$. If $(U_n)$ is a
decreasing basis of zero neighbourhoods in $E$, then for every $m$ there is
\[
  x_m \in co Z_m \cap (x+U_m)\, .
\]
In particular, $x_m\lora x$ in $E$, and therefore $x_n\lora x$ in some
$E_C$. Since $co\,Z_m\supset co\, Z_{m+1}$, the sequence $(x_n)_{n\ge m}$
is
contained in $co\, Z_m$ for all $m$, hence $x\in \ov{ co\, Z_m}^{E_C}$ for
all $m$.
\eop

As announced, we finally give a list of natural examples:

\bsk
\begin{exams}\label{exams2.8}\hfill

\bn
{\rm
1. A K\"othe matrix $A=(a_n)_n$  is a
sequence of scalar sequences satisfying $0<a_n(i) \le a_{n+1}(i)$ for each
$n$, $i$. The K\"othe echelon space associated with $A$ is the \Fre \
space defined by
\[
  \lambda_1(A) := \Blk x\in \R^{\N} \bmv p_n(x) := \ma{\tsum}\li {i}
   a_n(i) |x(i)| \text{ for all } n \Brk\, ,
\]
endowed with the metrizable, locally convex topology generated by the
fundamental system of seminorms $(p_n)$. The space $\lambda_1(A)$ is
nuclear if and only if for each $n$ there is $m>n$ such that $\bl a_n(i)/
a_m(i)\br_i\in \ell_1$. If $\lambda_1(A)$ is nuclear, then its strong dual
is the complete nuclear (DF)-space
\begin{align*}
  \lambda_1(A)'_b & = k_{\infty}(A)\\
                  & = \Blk u\in \R^{\N} \Bmv\ex n : \| u \|_n :=
                      \sup_{i} |u(i)| /a_n(i) < \infty \Brk \, .
\end{align*}

In fact, $\lambda_1(A)'_b = \rind_n \ell_\infty (1/a_n)$ is an (LB)-space
and its topology can be described by the seminorms
\[
  p_v(u) := \sup_{i}v(i) |u(i)|
\]
as $v$ varies in the set of nonnegative sequences such that $\bl a_n(i)
v(i)\br_i\in \ell_\infty$ for all $n$. For more details and
examples we refer to \cite{Bi} and \cite{MV}. One of the  most important
examples
is the space $s$ of rapidly decreasing sequences, which corresponds to
$\lambda_1(A)$ for the K\"othe matrix $(a_n(i)) := (i^n)$. The space
of Schwartz tempered distributions $S'$ on $\R^N$, is isomorphic
to $s'$. Moreover, $C^\infty[0,1]$ and $\cd(K)$, $K$ a compact subset of
$\R^N$ with nonempty interior, are isomorphic to $s$.

\bn
2. Let $K$ be a nonvoid compact subset of $\C^N$. The space
$H(K)$ of germs of holomorphic functions on $K$ is a complete nuclear
(DF)-space. Its topology is described as a countable inductive limit of
Banach spaces in the following way: let $(U_n)$ be a decreasing basis of
open neighbourhoods of $K$ satisfying $U_n\ssp \ov{U_{n+1}}$. Let
\[
  i_n: H^\infty (U_n)\lora H^\infty (U_{n+1})
\]
be the restriction map which is absolutely summing. Then $H(K) = \rind_n\,
H^\infty(U_n)$.

\bn
3. An (LF)-space $E = \rind_n\, E_n$ is a Hausdorff countable, inductive
limit of an increasing sequence of \Fre \ spaces. An (LF)-space is called
strict if every space $E_n$ is a closed topological subspace of
$E_{n+1}$. In this case, each $E_n$ is a closed topological subspace of
$E$, and every convergent sequence in $E$ is contained and converges in a
step $E_n$. Accordingly, if $E = \rind_n$ $E_n$ is a strict $(LF)$-space
and each step $E_n$ is nuclear, we can apply \Ban 's original theorem to
conclude that the domain of sums of each convergent series in $E$ is given
by the same formula. The most important example is the space of test
functions for distributions $\cd(\Om)$ on an open set $\Om\subset \R^N$.
By a result of Valdivia \cite{Va} and Vogt \cite{V}, $\cd(\Om)$
is in fact isomorphic to a countable direct sum $s^{(\N)}$ of copies of the
space $s$. We observe that the space $\cd_{(\omega)}(\Om)$ of test
functions
for ultradistributions of Beurling type on an open subset $\Om$ of $\R^N$
is also a strict (LF)-space (see \cite{BrMT}).

\bn
4. By the sequence space representation of Valdivia \cite{Va} and of Vogt
\cite{V},
the space of distributions $\cd'(\Om)$ on an open subset $\Om$ of $\R^N$
is isomorphic to a countable product $(s')^{\N}$ of copies of the
space $s'$. Therefore, it is a complete co-nuclear space. Moreover, by
the permanence properties of the strict Mackey condition (see
\cite[5.1.31]{BPC}), every convergent sequence in $\cd'(\Om)$ converges
locally, and our Theorem \ref{theo2.6} can be applied to $\cd'(\Om)$.
More spaces have a similar structure: the space of ultradistributions
$\cd'_{(\omega)}(\Om)$ of Beurling type and the space
$\ce_{\{\omega\}}(\Om)$ of ultradifferentiable functions of Roumieu type on
an open subset $\Om$ of $\R^N$ are also complete co-nuclear
spaces
which are the countable, projective limit of complete nuclear (DF)-spaces.
In particular, they are subspaces of countable products of complete nuclear
(DF)-spaces, and again the permanence properties of the strict Mackey
condition ensure that Theorem
\ref{theo2.6} can be applied. We refer to \cite{BrMT}.

\bn
5. The space $A(\Om)$ of real-analytic functions on an open subset $\Om$ of
$\R^N$ is endowed with a locally convex topology as follows:
$A(\Om) = \proj_K H(K)$, as $K$ runs over all the compact subsets of $\Om$
and $H(K)$ is defined as in 2. Deep results of Martineau from the 60's (cf.
\cite{BDom}) imply that $A(\Om)$ is an ultrabornological, countable,
projective limit of complete nuclear (DF)-spaces, and its dual $A(\Om)'_b$
is a complete nuclear
(LF)-space in which every convergent sequence is contained and converges in
a step. Therefore, $A(\Om)$ and $A(\Om)'_b$ are both complete co-nuclear
spaces to which Theorem \ref{theo2.6} can be applied
(again these spaces satisfy the strict Mackey condition by the argument
given in 4.).
}
\end{exams}

%3
\section{The converse}

As already mentioned, \Ban \ in \cite{B2} proved the converse of his
rearrangement theorem from \cite{B1}. Modifying his cycle of ideas, we will
obtain the following converse of our rearrangement theorem for nuclear
(DF)-spaces from section 2 (Theorem \ref{theo2.1}(a)).

\begin{theo}\label{theo3.1}
Let $E$ be a (DF)-space such that each convergent sequence converges
locally. If for each convergent series $\suk$ in $E$
\[
  S\bl\suk\br = \skinf u_k + \Gamma^{\bot}_{\loc} \bl\suk\br\, ,
\]
then $E$ is nuclear.
\end{theo}

We will first prove an appropriate characterization of nuclear (DF)-spaces
which seems to be interesting in its own right. If $X$ is an
$n$-dimensional Banach space, $\lambda_n$ the Lebesgue measure on $\R^n$
and $\phi: \R^n \lora X$ a linear bijection, then
\[
  \vol_n(A):= \lambda_n\bl \phi^{-1} A), \; A \text{ a Borel set in } X
\]
gives a measure on the Borel sets of $X$. Although $\vol_n(A)$ changes
with $\phi$, the ratio $\vol_n(A)$ $\vol_n(B)^{-1}$
for two Borel sets $A$ and $B$ is certainly
independent. For a linear bijection $T: X\lora Y$ between two
$n$-dimensional Banach spaces define
\[
  \vn (T) := \bl\vol_n(T B_X) / \vol_n(B_Y)\br^{1/n}\, ;
\]
for a linear operator $T: X\lora Y$ between two arbitrarily
normed spaces set
\[
  \vn (T) := \sup \vn (T: M\lora TM)\, ,
\]
where $n\in \N$ and the $\sup$ is taken with respect to all subspaces $M$
of $X$ such
that $\dim M = \dim TM = n$; put $\vn (T) = 0$ whenever rank $T< n$. For
this definition and the following (later very important) properties see
\cite[Lemma 1]{BB} and also \cite{Bl} :
\begin{itemize}
\item[(1)]
 $\vn (T) \le \| T\|$
\item[(2)]
 $\vn (TS) \le \vn (T)\, \vn (S)$
\item[(3)]
 $\vn (RTS) \le \| R\| \vn (T)\| S\|$
\item[(4)]
 $\vn (T) = \Bl \ma{\prod}\li{i=1}^n \delta_i(T)\Br^{1/n}$,  whenever
             $T$ acts between Hilbert spaces and $\delta_i$ stands for the
             $i$-th Kolmogorov number.
\item[(5)]
 $\vn (T) \ge h_n(T)$, where $h_n(T)$ stands for the $n$-th Hilbert number.
\end{itemize}

Moreover, denote for each $0<\varepsilon<1$ the class of all operators $T$
between normed spaces such that
\[
  \bl n^\varepsilon \vn (T) \br \in \e_\infty
\]
by $\fv_\varepsilon$. Part (b) of the next proposition is implicit in
\cite[Lemma 2]{BB}.

\bsk
\begin{prop}\label{prop3.2}
Let $0<\varepsilon<1$.
\begin{itemize}
\item[{\rm (a)}]
 Each composition of three {\rm 2}-summing operators belongs to
  $\fv_\varepsilon$.
\item[{\rm (b)}]
 The composition of $k\ge 5/\varepsilon$ operators in $\fv_\varepsilon$ is
nuclear.
\end{itemize}
\end{prop}

\Proof
(a) From Pietsch's factorization theorem we know that such a composition
factorizes through a nuclear operator $T: H_1\lora H_2$, $H_1$ and $H_2$
Hilbert spaces. Hence it suffices to check that $T$ belongs to
$\fv_\varepsilon$:
it is well-known that $\bl \delta_k(T)\br \in\e_1$, hence the assertion
is an immediate consequence of
\[
  n^\varepsilon v_n(T) \stackrel{(4)}{=} n^\varepsilon\Bl\prod_1^n
  \delta_i(T)\Br^{1/n}\le n^\varepsilon \frac{1}{n} \sum_{i=1}^n
  \delta_i(T) \, .
\]

\bn
(b) For every composition $T=T_k\, \circ\ldots\circ \, T_1$ of $k\ge
5/\varepsilon$ operators in $\fv_\varepsilon$ we have
\[
   \sup_{\e}\e^5 v_{\e}(T) \stackrel{(2)}{\le}\, \sup_{\e}\,\prod_{i=1}^k
   \e^\varepsilon v_{\e}(T_i) =: c<\infty\, .
\]
On the other hand, for each $n$
\begin{align*}
  \delta_n(T) & \le \Bl \prod^n_{\e=1}\delta_{\e}(T)\Br^{1/n}\\[1ex]
         & \le n\Bl \prod^n_{\e=1}h_{\e}(T)\Br^{1/n}\q \text{(see the proof
             of \cite[12.12.3]{P1})}\\[1ex]
         & \stackrel{(5)}{\le}\Bl \prod^n_{\e=1}n v_{\e}(T)\Br^{1/n}\, .
\end{align*}
Together we obtain
\[
   \delta_n(T) n^3  \le \Bl \prod^n_{\e=1}\, \frac{n^4}{\e^5}\, \e^5
                v_{\e}(T)\Br^{1/n}
                \le c \Bl \prod^n_{\e=1}\, \frac{n^4}{\e^5}\Br^{1/n}\lora
                0, \; n\lora \infty\, ,
\]
which shows that $\bl \delta_n(T) n\br \in \e_1$, hence by
\cite[11.12.2]{P2}
the sequence of approximation numbers $\bl a_n(T)\br\in \e_1$. It is
well-known that this assures the nuclearity of $T$ (see
\cite[18.6.3]{P2}).
\eop

\bsk
As an immediate consequence we obtain the following corollary.
\begin{cor}\label{cor3.3}
Let $E$ be a locally convex space. Then the following are equivalent:
\begin{itemize}
\item[{\rm(1)}]
  $E$ is nuclear.
\item[{\rm(2)}]
  $\fa\, 0 < \varepsilon < 1\, (\ex\; 0 < \varepsilon < 1) \, \fa\, p\in
   cs\, (E)\; \ex\; p\le q\in cs\, (E) : \pi^p_q\in \fv_\varepsilon$.
\end{itemize}

\noi
Dually, the following equivalence holds:
\begin{itemize}
\item[{\rm(1')}]
 $E$ is co-nuclear.
\item[{\rm(2')}]
 $ \fa\, 0 < \varepsilon < 1 \,(\ex\; 0< \varepsilon < 1) \, \fa\, B\in
 \cb(E)\; \ex\; B\subset C\in \cb (E): i^C_B\in \fv_\varepsilon$.
\end{itemize}
\end{cor}

\bn
The implication (2) $\car$ (1) was stated in \cite[Lemma 2]{BB}. For our
purposes we need a characterization of nuclear (DF)-spaces via volume
numbers which uses bounded sets and continuous seminorms simultaneously --
based on ideas of \cite[6.4.2]{Ju} we prove:

\begin{prop}\label{prop3.4}
For every (DF)-space $E$ the following are equivalent:
\begin{itemize}
 \item[{\rm(1)}]
$E$ is nuclear.
\item[{\rm(2)}]
 $\fa\, 0<\varepsilon<1 \; (\ex\; 0<\varepsilon<1)\, \fa\, B\in\cb (E)\,
 \fa\,  p\in cs\,(E): \pi_p i_B\in \fv_\varepsilon$.
\end{itemize}
\end{prop}

\Proof
The proof of (1) $\car$ (2) is easy: clearly, in every nuclear space $E$
each
mapping $\pi_p i_B$ can be written as a composition of
three nuclear maps,
hence $\pi_p i_B\in \fv_\varepsilon$ for each $0<\varepsilon<1$ (Lemma
3.2). Conversely, assume that (2) holds for some $0<\varepsilon<1$
and take a
set $B\in \cb(E)$. By Corollary \ref{cor3.3} it suffices to show that
$i^C_B\in \fv_\varepsilon$ for some $B\subset C\in \cb(E)$. Take a
fundamental system $(B_n)$ of bounded sets such that
\[
  B_1 := B \q\text{and}\q 2B_n\subset B_{n+1} \q\text{for all } n\, .
\]
Assume that
\[
  i_n := i^{B_n}_{B_1}\not\in \fv_\varepsilon \text{ for all }n>1\, .
\]
By definition, for every $n$ there is a $k_n$-dimensional subspace
$M_n$ of $E_B$ such that
\begin{align*}
 & \dim M_n = \dim i_nM_n = k_n>k_{n-1}\\
 & k^\varepsilon_n v_{k_n}\bl i_n: M_n\lora i_n M_n\br
   \ge n\, .
\end{align*}
It is well-known (see e.g. \cite[6.3]{DF}) that for each $n$ there is a
subspace $G_n$ of some $\e^{m(n)}_\infty$ and a linear bijection
\begin{align*}
 &  j_n        : i_nM_n\lora G_n\\
 &  \| j_n \|  = 1\\
 &  \|j_n^{-1}  : G_n\lora i_nM_n\| \le 2\, .
\end{align*}
By the Hahn-Banach theorem $j_n$ has an extension $\wt{j}_n : E_{B_n}\lora
\e^{m(n)}_\infty$ with equal norm, and we know fom \cite[4.3.12]{Ju} that
there is an operator $Q_n\in \cl \bl E$, $\e_\infty^{m(n)}\br$ such that
\begin{align*}
  &  \| Q_{n \mid E_{B_n}}\| \le 2 \\[1ex]
  &  \| Q_n h -\wt{j}_n h\|_\infty \le 1/4 \| h\|_{E_{B_n}} \text{ for all }
     h\in i_nM_n\, ;
\end{align*}
summarizing, we get the following diagram:

\begin{center}
\special{em:linewidth 0.4pt}
\unitlength 1.00mm
\linethickness{0.4pt}
\begin{picture}(55.33,44.00)
\put(25.21,12.45){\oval(1.56,2.22)[b]}
%\vector(25.99,12.23)(25.99,18.90)
\put(25.99,18.90){\vector(0,1){0.2}}
\emline{25.99}{12.23}{1}{25.99}{18.90}{2}
%\end
\put(4.21,12.45){\oval(1.56,2.22)[b]}
%\vector(4.99,12.23)(4.99,18.90)
\put(4.99,18.90){\vector(0,1){0.2}}
\emline{4.99}{12.23}{3}{4.99}{18.90}{4}
%\end
\put(25.55,32.45){\oval(1.56,2.22)[b]}
%\vector(26.33,32.23)(26.33,38.90)
\put(26.33,38.90){\vector(0,1){0.2}}
\emline{26.33}{32.23}{5}{26.33}{38.90}{6}
%\end
\put(50.88,12.44){\oval(1.56,2.22)[b]}
%\vector(51.66,12.22)(51.66,18.89)
\put(51.66,18.89){\vector(0,1){0.2}}
\emline{51.66}{12.22}{7}{51.66}{18.89}{8}
%\end
\put(5.00,25.00){\makebox(0,0)[cc]{$E_B$}}
\put(26.67,24.67){\makebox(0,0)[cc]{$E_{B_n}$}}
\put(51.66,25.33){\makebox(0,0)[cc]{$\e^{m(n)}_\infty$}}
\put(5.00,5.67){\makebox(0,0)[cc]{$M_n$}}
\put(26.00,6.00){\makebox(0,0)[cc]{$i_nM_n$}}
\put(51.66,6.33){\makebox(0,0)[cc]{$G_n$}}
\put(12.12,6.79){\oval(2.22,1.56)[l]}
%\vector(11.90,6.01)(18.57,6.01)
\put(18.57,6.01){\vector(1,0){0.2}}
\emline{11.90}{6.01}{9}{18.57}{6.01}{10}
%\end
\put(12.45,25.46){\oval(2.22,1.56)[l]}
%\vector(12.23,24.68)(18.90,24.68)
\put(18.90,24.68){\vector(1,0){0.2}}
\emline{12.23}{24.68}{11}{18.90}{24.68}{12}
%\end
\put(14.67,7.34){\makebox(0,0)[cb]{$\scriptstyle{i_n}$}}
\put(14.99,26.00){\makebox(0,0)[cb]{$\scriptstyle{i_n}$}}
\put(26.33,44.00){\makebox(0,0)[cc]{$E$}}
\put(40.67,36.33){\makebox(0,0)[cb]{$\scriptstyle{Q_n}$}}
\put(38.34,26.32){\makebox(0,0)[cb]{$\scriptstyle{\wt{j}_n}$}}
\put(37.67,7.33){\makebox(0,0)[cb]{$\scriptstyle{j_n}$}}
%\vector(32.67,6.00)(44.67,6.00)
\put(44.67,6.00){\vector(1,0){0.2}}
\emline{32.67}{6.00}{13}{44.67}{6.00}{14}
%\end
%\vector(32.67,25.00)(44.67,25.00)
\put(44.67,25.00){\vector(1,0){0.2}}
\emline{32.67}{25.00}{15}{44.67}{25.00}{16}
%\end
%\vector(33.00,39.00)(44.67,31.67)
\put(44.67,31.67){\vector(3,-2){0.2}}
\emline{33.00}{39.00}{17}{44.67}{31.67}{18}
%\end
\put(55.33,6.33){\makebox(0,0)[cc]{.}}
\end{picture}
\end{center}

\sn
We prove that $Q_n$ is invertible on $i_nM_n$ and
\[
  \| Q^{-1}_n: Q_n i_nM_n \lora i_nM_n\| \le 4\, ;
\]
indeed, for $h\in i_nM_n$
\begin{align*}
   \| h\|_{E_{B_n}} & = \| j^{-1}_n j_nh\|_{E_{B_n}} \le 2 \|
                        j_nh\|_\infty\\[1ex]
                    & \le 2 \|j_nh - Q_nh\|_\infty + 2 \| Q_nh\|_\infty \le
                       1/2 \|h \|_{E_{B_n}} + 2\| Q_nh \|_\infty \, ,
\end{align*}
hence
\[
  \| h \|_{E_{B_n}} \le 4 \| Q_nh \|_\infty \, .
\]
The norm estimate for $Q^{-1}_n$ can now be used to show that for each $n$
\begin{align*}
  v_{k_n}\bl M_n \stackrel{i_n}{\lora}\, i_nM_n \br
 &    =   v_{k_n}\bl M_n\stackrel{i_n}{\lora}\, i_nM_n
       \stackrel{j_n}{\lora}\, G_n\stackrel{j^{-1}_n}{\lora}\,
        i_nM_n\br\\[0.5ex]
 &   \le 2 v_{k_n}\bl M_n\stackrel{i_n}{\lora}\, i_nM_n
       \stackrel{Q_n}{\lora}\, Q_ni_nM_n
       \stackrel{Q^{-1}_n}{\lora}\,i_nM_n
       \stackrel{j_n}{\lora}\,G_n\br\\[0.5ex]
 &   \le  8 v_{k_n}\bl M_n\stackrel{i_n}{\lora}\, i_nM_n
       \stackrel{Q_n}{\lora}\, Q_ni_nM_n \br \, ,
\end{align*}
hence
\begin{align*}
  1/8n & \le 1/8 k^\varepsilon_n v_{k_n}\bl M_n\stackrel{i_n}{\lora}\,
            i_nM_n\br\\
       & \le k_n^\varepsilon v_{k_n}\bl
            M_n\stackrel{i_n}{\lora}\,i_nM_n
             \stackrel{Q_n}{\lora}\, Q_ni_nM_n\br\, ,
\end{align*}
and finally
\[
   1/8 n \le k^\varepsilon_n v_{k_n}\bl E_B\stackrel{i_B}{\hra}\,
    E\stackrel{Q_n}{\lora}\, \e^{m(n)}_\infty\br\, .
\]
Define
\[
  V := \bigcap_{\N} Q^{-1}_n B_{\e^{m(n)}_\infty}\subset E\, ;
\]
$V$ absorbs bounded sets and is hence a zero neighbourhood in the
(DF)-space $E$: clearly, for each $k$ there is $\lambda>0$ such that
\[
  B_k \subset \lambda \bigcap_{n=1}^{k} Q^{-1}_n B_{\e^{m(n)}_\infty}\, ,
\]
and for $n>k$ we have $2 B_k \subset B_n$, hence
\[
  2 Q_n B_k \subset Q_n B_n \subset 2B_{\e^{m(n)}_\infty}\, .
\]
Now observe that $Q_nV \subset B_{\e^{m(n)}_\infty}$ for each $n$ which assures
that there are operators\linebreak[4] $\wh{Q}_n: E_V\lora \e^{m(n)}_\infty$
with
$\|\wh{Q}_n \| \le 1$ and such that the following diagram commutes:

\bsk
\begin{center}
\special{em:linewidth 0.4pt}
\unitlength 1.00mm
\linethickness{0.4pt}
\begin{picture}(51.99,26.99)
\put(5.00,25.67){\makebox(0,0)[cc]{$E_B$}}
\put(25.00,25.34){\makebox(0,0)[cc]{$E$}}
\put(51.99,26.00){\makebox(0,0)[cc]{$\e^{m(n)}_\infty$}}
\put(12.45,26.13){\oval(2.22,1.56)[l]}
%\vector(12.23,25.35)(18.90,25.35)
\put(18.90,25.35){\vector(1,0){0.2}}
\emline{12.23}{25.35}{1}{18.90}{25.35}{2}
%\end
\put(14.99,26.67){\makebox(0,0)[cb]{$\scriptstyle{i_B}$}}
\put(36.67,26.99){\makebox(0,0)[cb]{$\scriptstyle{Q_n}$}}
%\vector(31.00,25.67)(43.00,25.67)
\put(43.00,25.67){\vector(1,0){0.2}}
\emline{31.00}{25.67}{3}{43.00}{25.67}{4}
%\end
%\vector(24.83,20.67)(24.83,8.67)
\put(24.83,8.67){\vector(0,-1){0.2}}
\emline{24.83}{20.67}{5}{24.83}{8.67}{6}
%\end
\put(26.66,15.00){\makebox(0,0)[lc]{$\scriptstyle{\pi_V}$}}
\put(24.83,3.67){\makebox(0,0)[cc]{$E_V$}}
%\vector(31.00,8.84)(43.00,19.50)
\put(43.00,19.50){\vector(1,1){0.2}}
\emline{31.00}{8.84}{7}{43.00}{19.50}{8}
%\end
\put(39.34,13.67){\makebox(0,0)[lc]{$\scriptstyle{\wh{Q}_n}$}}
\end{picture}
\end{center}

\vspace*{-1ex}
\noi
(as usual, $E_V = E_{m_V}$ and $\pi_V = \pi_{m_V}$).
Since by assumption $\bl m^\varepsilon v_m(\pi_V i_B)\br\in \e_\infty$,
we obtain a contradiction:
\begin{align*}
  1/8 n & \le k^\varepsilon_n v_{k_n} (Q_ni_B)\\
        & \le k^\varepsilon_n v_{k_n} \bl \wh{Q}_n \pi_V i_B\br \le
           k^\varepsilon_n v_{k_n}(\pi_V i_B)\, .
\end{align*}
\eop

\ssk
We need the following deep result of \Ban \ \cite[Lemma 4]{B2}
(which among others is based on Milman's quotient subspace theorem from
local Banach space theory). Here $gp(A)$ denotes the subgroup generated by
a subset $A$ of a group $G$; moreover, for a seminorm $p$ on a vector space
$E$ we write $B_p$ for the closed unit ball with respect to $p$, and
$d_p(x,H)$
for the distance of $x\in E$ and $H\subset E$ with respect to $p$.

\begin{theo}\label{theo3.5}
Let $E$ be a vector space and $q\ge p$ two seminorms on $E$ such that
\[
   \bl \pi^p_q: E_q\lora E_p\br \not\in \fv_{0.1}\, .
\]
Then for each finitely generated subgroup $G\subset E$, each $a\in \spa G$
and each $\gamma > 0$ there is a finitely generated group $G_1
\ssp G$ satisfying
\begin{align*}
(1)   & \q  G_1 = gp (G_1\cap B_q)\, ,\\
(2)   & \q  d_p(a, G_1) \ge d_p(a,G)-\gamma\, .
\end{align*}
\end{theo}

\sn
The following consequence will be crucial.

\begin{cor}\label{cor3.6}
Let $E$ be a locally convex space, $B\in\cb(E)$ and $p\in cs (E)$ such that
$\pi_p i_B\not\in \fv_{0.1}$. Then there is a subgroup $G $ of $E$ such
that
\begin{align*}
(1)   & \q  G = gp (G\cap 1/m\, B) \text{ for all } m\, ,
\\
(2)   & \q  1/2\, G \not\subset \ov{G}^p\, .
\end{align*}
\end{cor}

\Proof
Take $a\in E_B$ such that $p(a)>0$ $(\pi_p i_B\ne 0!)$. Without loss of
generality we may assume that
\[
  p(a) = 2 \q\text{and}\q p\le \|\cdot \|_B \text{ on } E_B\, .
\]
Define the finitely generated subgroup $G_0 := 2a\Z$ of $E_B$. Obviously,
\[
  p(a-2a\, z) = 2|1-2z|\text{ for all } z\in \Z\, ,
\]
hence $d_p(a,G_0) = 2$. Clearly, $(E_B)_p$ is an isometric subspace of
$E_p$, hence the canonical map
\[
  \bl \pi^p_B: E_B\lora (E_B)_p\br \not\in \fv_{0.1}\, .
\]
By \Ban 's Theorem 3.5 there is a finitely generated subgroup $G_0\subset
G_1 \subset E_B$ such that
\begin{align*}
  &   G_1 = gp\, (G_1\cap B)\, ,\\
  &   d_p(a,G_1) \ge d_p (a,G_0) - 1/2\, .
\end{align*}
Next apply \Ban 's result 3.5 to $G_1$, $2\|\cdot \|_B$, $2p$ and $1/2$:
there is a finitely generated subgroup $G_2\supset G_1$ such that
\begin{align*}
  &   G_2 = gp \Bl G_2\cap \frac{1}{2}B\Br\, , \\[1ex]
  &   d_{2p}(a,G_2) \ge d_{2p} (a,G_1) - \frac{1}{2}\, , \text { hence }
      d_p(a,G_2)\ge d_p(a,G_1) - \frac{1}{2^2}\, .
\end{align*}
Proceeding this way, one gets an increasing sequence of finitely generated
subgroups\linebreak[4] $G_n\supset G_{n-1}$ such that
\begin{align*}
  G_n        &   =  gp \Bl G_n\cap \frac{1}{n} B\Br\, ,\\[1ex]
  d_p(a,G_n) &  \ge d_p (a,G_{n-1}) - \frac{1}{2^n}\, .
\end{align*}
Define the subgroup
\[
  G:= \bigcup_n G_n
\]
of $E$. Then for all $m$
\[
  G = gp\, \bl G\cap 1/m\, B\br\, ;
\]
indeed, for $m$ and $x\in G$ we have $x\in G_n$ for some $n>m$, hence
\[
  x=\sum_{\text{finite}} g_i\text{ with } g_i\in G_n\cap 1/n\, B \subset G
    \cap 1/m\,  B.
\]
Moreover, $2a\in G_0 \subset G$, but
\begin{align*}
  d_p(a,G) &  =   \inf \blk d_p (a,G_k) \bmv k\in\N \brk\\
           &  \ge \inf \Blk d_p(a,G_0) - \sum_{\e =1}^k\,
              \frac{1}{2^{\e}} \Bmv k\in\N \Brk\\[1ex]
           &  = 2-\sum_{\e =1}^\infty \, \frac{1}{2^{\e}} = 1\, ,
\end{align*}
hence $a\not\in \ov{G}^p$.
\eop

\ssk
Finally, we are prepared to give a {\sc proof} of Theorem \ref{theo3.1}:
Assume that $E$ is not nuclear. Then we know from Proposition \ref{prop3.4}
that
\[
  \pi_p i_B\not\in \fv_{0.1}
\]
for some $B\in \cb (E)$ and $p\in cs\, (E)$. Hence by Corollary \ref{cor3.6}
there are a subgroup $G$ of $E$ and an $a\in E$ such that
\begin{itemize}
\item[(1)]
 $G = gp\, (G\cap 1/m\, B)$ for all $m$,
\item[(2)]
 $2a\in G$, but $a\not\in \ov{G}^p$.
\end{itemize}
Take a fundamental system $(B_n)$ of $\cb(E)$ such that $B_1=B$. By (1) we
have for each $m$ the following finite representations of $2a$:
\[
  2a = \sum_{i=1}^{s(m)}\, w_i^m\, , \, w_i^m\in G \cap 1/m\, B\, .
\]
Define the series $\tsum(u_k)$ by
\[
   (u_k) := \bl w_1^1, -w_1^1, \ldots, w_{s(1)}^1, - w_{s(1)}^1, w_1^2, -
    w_1^2, \ldots, w_{s(2)}^2, - w_{s(2)}^2, \ldots\br \, .
\]
Obviously, the series $\suk$ is convergent in $E_{B_1}$ and we have
(convergence in $E_{B_1}$)
\[
  0 = \sum_{k=1}^\infty \, u_k\, ,
\]
and
\[
  2a\in A\subset \bigcup_n A_{B_n}\, .
\]
We now show that $a\not\in \ma{\bigcup}\li n A_{B_n}$ which then
contradicts the fact that (by assumption and Remark \ref{rem2.4})
\[
  S\bl \suk \br = \bigcup_n A_{B_n} = \Gamma^{\bot}_{\loc}
\]
is convex: assume that $a\in A_{B_n}$ for some $n$. Since for all $m$
\[
  Z_m \subset gp\, \blk w_i^k \bmv i= 1,\ldots, s(k) \q \text{and}\q
  k\in\N \brk \subset G\, ,
\]
we get that
\[
  a\in A_{B_n} = \bigcap_m \ov{Z}_m^{E_{B_n}} \subset
\ov{G}^{E_{B_n}}\subset \ov{G}^p\, ,
\]
a contradiction.\eop

\bsk
An easy analysis of the preceding proof gives slightly more:

\begin{rem}\label{rem3.7}
Let $E$ be a locally convex space such that for some $B\in \cb(E)$ and
$p\in cs\, (E)$
\[
  \pi_p i_B\not\in \fv_{0.1}\, .
\]
Then there is some $E_B$-convergent series $\suk$ and there is $a\in E_B$
such that
\begin{align*}
  2a\in \bigcap_m Z_m & = A\\
  a\not\in \bigcap_m\ov{Z}_m^p & =: A_p\, .
\end{align*}
\end{rem}
This remark shows that our method also gives a new proof of the converse of
\Ban 's rearrangement theorem for metrizable spaces; indeed, if $E$ is a
nonnuclear, metrizable space such that for each convergent series $\suk$
\[
   S\bl \suk\br = \sum_{k=1}^\infty u_k + \Gamma^{\bot}\bl\suk\br\, ,
\]
then by Lemma \ref{lemma1.2} for each such series
\[
   A_E\bl\suk\br
\]
is a subspace. But this contradicts Remark \ref{rem3.7},
since $A_p\ssp A_E\ssp A$ and the following counterpart of Proposition
\ref{prop3.4} assures the existence of an appropriate $B$ and
$p$.

\begin{prop}\label{prop3.8}
For every metrizable space the statements {\rm (1)} and {\rm (2)} of
Proposition {\rm \ref{prop3.4}} are equivalent.
\end{prop}

\Proof
That (1) implies (2) follows as in Proposition \ref{prop3.4}. Hence assume
that $E$ satisfies (2) for some $0< \varepsilon< 1$. Let $p\in cs\, (E)$; by
Corollary \ref{cor3.3} it suffices to show that there is some $p\le q\in
cs (E)$ with $\pi_q^p\in \fv_\varepsilon$. Fix a fundamental system $(p_n)$
of $cs\, (E)$ such that
\[
  p_n\ge p_1 = p \text{ for all } n\, ,
\]
and assume that
\[
  \pi_n := \pi_{p_n}^p \not\in \fv_\varepsilon \text{ for all } n>1\, .
\]
By definition for every $n$ there is a $k_n$-dimensional subspace
$M_n$ of $E_{p_n}$ such that
\begin{align*}
   &    \dim M_n = \dim \pi_n M_n = k_n> k_{n-1}\\[0.5ex]
   &    k_n^{\varepsilon} v_{k_n} (\pi_n: M_n\lora \pi_n M_n) \ge
        n\, .
\end{align*}
The canonical embedding $i_n: M_n\hra E_{p_n}$ being of finite rank
has a finite representation
\[
   i_n = \tsum \varphi'_j \otimes \pi_{p_j}(x_j)\in M'_n \otimes
   E_{p_n}\, .
\]
Define the operator
\[
   R_n := \tsum \varphi'_j\otimes x_j: M_n\lora E\, ;
\]
clearly, $i_n = \pi_{p_n}R_n$ for all $n$, hence
\begin{align*}
  &   k_n^\varepsilon v_{k_n}(\pi_pR_n) = k_n^\varepsilon v_{k_n}(\pi_n
   \pi_{p_n}R_n)\\[0.5ex]
  &   = k_n^\varepsilon v_{k_n}(\pi_n i_n) = k_n^\varepsilon v_{k_n}(\pi_n:
   M_n\lora \pi_nM_n)\ge n\, .
\end{align*}
Then the closed, absolutely convex hull $B$ of $\ma{\bigcup}\li n
R_nB_{M_n}$ is a bounded set in $E$; indeed, for $k$ given, the set
$\{ x\in E \mid p_k(x) \le 1\}$ absorbs the bounded set $\ma{\bigcup}\li
1^k R_n B_{M_n}$, and for $n>k$
\[
  R_n B_{M_n}\subset \blk x\in E\bmv p_n(x) \le 1\brk \subset \blk
  x\in E\bmv p_k(x) \le 1\brk\, .
\]
Observe now that $R_nB_{M_n} \subset B$ for all $n$, hence there are
operators $\wh{R}_n : M_n \lora E_B$ with norm $\le 1$ such that

\bsk
\begin{center}
\special{em:linewidth 0.4pt}
\unitlength 1.00mm
\linethickness{0.4pt}
\begin{picture}(54.77,23.65)
\put(29.88,9.78){\oval(1.56,2.22)[b]}
%\vector(30.66,9.56)(30.66,16.23)
\put(30.66,16.23){\vector(0,1){0.2}}
\emline{30.66}{9.56}{1}{30.66}{16.23}{2}
%\end
\put(7.34,22.45){\makebox(0,0)[cc]{$M_n$}}
\put(31.12,22.56){\makebox(0,0)[cc]{$E$}}
\put(54.77,22.44){\makebox(0,0)[cc]{$E_p$}}
\put(43.01,23.65){\makebox(0,0)[cb]{$\scriptstyle{\pi_p}$}}
%\vector(37.34,22.33)(49.34,22.33)
\put(49.34,22.33){\vector(1,0){0.2}}
\emline{37.34}{22.33}{3}{49.34}{22.33}{4}
%\end
%\vector(12.67,16.06)(24.34,8.73)
\put(24.34,8.73){\vector(3,-2){0.2}}
\emline{12.67}{16.06}{5}{24.34}{8.73}{6}
%\end
\put(18.29,23.16){\makebox(0,0)[cb]{$\scriptstyle{R_n}$}}
%\vector(12.62,21.84)(24.62,21.84)
\put(24.62,21.84){\vector(1,0){0.2}}
\emline{12.62}{21.84}{7}{24.62}{21.84}{8}
%\end
\put(30.45,3.78){\makebox(0,0)[cc]{$E_B$}}
\put(32.01,12.22){\makebox(0,0)[lc]{$\scriptstyle{i_B}$}}
\put(17.23,11.45){\makebox(0,0)[rc]{$\scriptstyle{\wh{R}_n}$}}
\end{picture}
\end{center}

\sn
commutes. But then
\[
   n\le k_n^\varepsilon v_{k_n} (\pi_p R_n) = k_n^\varepsilon v_{k_n}
   (\pi_p i_B \wh{R}_n) \le k_n^\varepsilon v_{k_n}(\pi_p i_B)
\]
contradicts the fact that $\pi_p i_B\in \fv_\varepsilon$.
\eop

%xx

\end{document}